\documentclass[12pt]{amsart}
\usepackage{amscd,amsmath,amssymb,amstext,amsfonts,amsthm,latexsym}
\usepackage{verbatim,multirow}
\usepackage[english]{babel}
\usepackage[latin1]{inputenc}

\usepackage{url}

\usepackage[a4paper, top=4cm, bottom=3.8cm, left=3cm, right=3cm]{geometry}

\usepackage{hyperref}

\usepackage{graphics}
\usepackage{color}

\newtheorem{theorem}{Theorem}[section]

\newtheorem{corollary}[theorem]{Corollary}
\newtheorem{proposition}[theorem]{Proposition}

\theoremstyle{remark}
\newtheorem{remark}[theorem]{Remark}

\theoremstyle{definition}
\newtheorem{definition}[theorem]{Definition}

 \thanks{}
																			
\title{ The fundamental group and torsion group of Beauville surfaces}
\makeatletter

 \def\@textbottom{\vskip \z@ \@plus 1pt}
 \let\@texttop\relax
\@addtoreset{equation}{section}
\makeatother
\author{Ingrid Bauer, Fabrizio Catanese, Davide Frapporti }

\date{\today}

\begin{document}
\begin{abstract}
We give a survey on the fundamental group of surfaces isogenous to a higher product. If the surfaces are regular, e.g. if they are Beauville surfaces, the first homology group is a finite group. We present a MAGMA script which calculates the first homology groups of regular surfaces isogenous to a product.
\end{abstract}

\maketitle

\section*{Introduction}
One aim of this note is to give an overview of what is known about the fundamental groups and more specifically about the 
first homology groups of Beauville surfaces and of their generalizations.  In particular, we wrote a MAGMA script which  calculates the torsion group of  Beauville surfaces, or more generally of regular surfaces isogenous to a higher product.

Fundamental groups of algebraic varieties are a very interesting area of research and they are still quite mysterious objects, since usually it may be very hard to determine them. On the other hand, if one can calculate them in some specific cases (e.g., for surfaces of general type) they are  a powerful tool to conclude that two surfaces (of general type) are not contained in the same connected component of their moduli space.

Here we will restrict ourselves to the following situation:
let $C_1$ and $C_2$ be projective algebraic curves  of respective genera $g_1, g_2$ at least 2 and 
let $G$ be a subgroup of the automorphism group $\mathrm{Aut}(C_1\times C_2)$. We denote by  $X$ 
the quotient $(C_1\times C_2)/G$; when  $X$ is  singular,
we denote by $S$ the minimal resolution of singularities  of $X$.

We define $G^0:=G \cap (\mathrm{Aut}(C_1)\times \mathrm{Aut}(C_2))$. Then (compare  \cite{Cat00} for this and the following assertions)
$G^0$ is a subgroup  of index at most $2$  inside  $G <\mathrm{Aut}(C_1\times C_2)$ 
and  acts on each factor and ``diagonally'' on the product $C_1\times C_2$ (i.e., for each $\gamma \in G^0$, we can write $\gamma= (\gamma_1, \gamma_2)$,
and $\gamma (x_1, x_2) = (\gamma_1 (x_1) , \gamma_2 (x_2))$). 

By \cite{Cat00}, it is always possible to 
assume that $G^0$ acts faithfully on both factors: in this case we say that
$(C_1\times C_2)/G$ is a \textit{minimal realization}, and this minimal realization is moreover unique. 
From now on we shall assume that we have indeed a minimal realization. 
Then (see  \cite{Cat00}) there are two cases. 

 The \textit{mixed case} is the case where  the 
action of $G$ exchanges the two factors: in this case $C_1\cong C_2$ and $G^0 \lhd_2 G$.
 
The \textit{unmixed case} is the case where $G=G^0$, and  its two projections into $\mathrm{Aut}(C_j),$ $ j=1,2,$ are injective.
Identifying  then $G$ to a subgroup $G \subset \mathrm{Aut}(C_j)$, for $j=1,2$,   we say that $G$ acts diagonally: this means that we view $G $ as the   diagonal
subgroup inside $ G \times G \subset (\mathrm{Aut}(C_1)\times \mathrm{Aut}(C_2))$ .\\
In the unmixed case, $S$ is called a \textit{product-quotient surface} and $X$ is its \textit{quotient model};
while in  the mixed case,  $S$ is a \textit{mixed surface} and $X$ is a \textit{mixed quotient}.

In this note we want to focus on the \textsl{unmixed case}, so from now on we implicitly assume that $G=G^0$.
Some of the results can be extended/generalized to the mixed case, but the approach is slightly different.
We refer to \cite{Fra11, FP13} for further details.

In the first section we describe product-quotient surfaces and the algebraic data which determine them. In the second section we comment on the results on the fundamental group of product-quotient surfaces and their higher-dimensional analogues. The third section is dedicated to the description of a MAGMA script which, given as input a finite group $G$ and two ordered tuples 
$$T_1=(n_1, \ldots , n_r),  \quad T_2=(m_1, \ldots , m_s),$$

\noindent has as output  

\begin{itemize}
\item 
one representative   surface $S = (C_1 \times C_2)/G$ for each irreducible connected component of the moduli space, where $G$ is acting freely with signatures $(T_1,T_2)$, 
\item
for each such representative, the first homology group $H_1(S, \mathbb{Z})$ of $S$.
\end{itemize}

Finally, in the last section,  we show  several concrete and explicit calculations, obtained  through a
direct application of  our program.

\section{Product-Quotient Surfaces}
Let $G$ be a finite group and let  $S$ be a \textit{product-quotient surface} with \textit{quotient model} $X =(C_1 \times C_2)/G$.
If $G$ acts freely, then $X$ is smooth, and we have the following:

\begin{definition}
A surface $S$ is said to be \textit{isogenous to a (higher) product}  if $S$ is the quotient $(C_1 \times C_2)/G$, where $g_i : = genus (C_i) \geq 2$, 
and $G$ is a finite group acting freely on $C_1 \times C_2$.
\end{definition}

The word ``higher'' means that the respective genera of $C_1, C_2$ are $\geq 2$, in particular this implies that $S$ is of general type
(and with ample canonical divisor).
However, for commodity, from now on we shall drop  the  word ``higher''.

In the last years, a huge amount of new surfaces of general type with $p_g=q$ have been constructed as the quotient 
of a product of two curves by the action of a finite group; see
\cite{BC04, BCG08, BCGP08, BP10, BP13, Fra11, FP13} for $p_g=0$, 
\cite{CP09, Pol07, Pol09, MP10, FP13} for $p_g=1$,
\cite{Penegini, zuc} for $p_g=2$.  

\noindent In particular, we have a complete classification of surfaces isogenous to a product with $p_g=q$. The case where $G$ does not act freely is still open (and particularly difficult in the regular case $q \neq 0$), in spite of several  results. A difficulty which is peculiar 
of  the regular case is that the following can happen: the minimal resolution $S$ of singularities  $X$ may not be minimal surface. Up to now there are almost no techniques to decide whether there are exceptional curves of the first kind (i.e., smooth rational curves with self intersection $(-1)$) on $S$ and how to find them explicitly (cf. \cite{BP13}).

\

The description of product- quotient surfaces  is accomplished through the theory 
  of Galois coverings between projective curves (also named ``Riemann surfaces'').

\begin{definition}\label{gv}
Let  $g\geq 0 \mbox{ and  } m_1, \ldots, m_r > 1$
be integers.
The 
\textit{orbifold surface group of signature $(g;m_1, \ldots, m_r)$} is defined as:
\begin{equation*}
\begin{array}{lr}
\mathbb{T}(g;m_1,\ldots ,m_r):=
\langle a_1,b_1,\ldots, a_{g},b_{g},c_1, \ldots, c_r \mid& \\ 
 \qquad \qquad  \qquad \qquad  \qquad \qquad \qquad  
 c_1^{m_1}, \ldots, c_r^{m_r},\prod_{i=1}^{g} [a_i,b_i]\cdot c_1 \cdots c_r\rangle\,.
\end{array}
\end{equation*}

 Given  a finite group $H$, a \textit{generating vector} for $H$ of signature $(g;m_1,\ldots ,m_r)$ is a $(2g+r)$-tuple of
elements of $H$:
$$V:=(d_1,e_1,\ldots, d_{g},e_{g};h_1, \ldots, h_r)$$
such that 
\begin{itemize}
\item $V$ generates $H$, 
\item $\prod_{i=1}^{g}[d_i,e_i]\cdot h_1\cdot h_2\cdots h_r=1$ and 
\item there 
exists a permutation $\sigma \in \mathfrak{S}_r$ such that $\mathrm{ord}(h_i)=m_{\sigma(i)}$ for $i=1,\ldots, r$.
\end{itemize}
If $g=0$, then $V:=(h_1,\ldots, h_r)$ is called a \textit{spherical system of generators} of $H$ of signature $(m_1,\ldots,m_r)$.
\end{definition}

To give a generating vector of signature $(g;m_1,\ldots, m_r)$ for a finite group $H$ is equivalent to
give an \textit{appropriate orbifold homomorphism} 
$$\psi \colon \mathbb T (g;m_1,\ldots, m_r)\longrightarrow H,$$ i.e.,
a surjective homomorphism $\psi$  such that $\psi(c_i)$ has order $m_i$.

\begin{remark}\label{RET}
By Riemann's existence theorem (see \cite{Survey}), any curve $C$ of genus $g$ together with an action of
a finite group $H$ on it, such that $C/H$ is a curve $C'$ of genus $g'$, is determined (modulo  automorphisms)
by the following data:
\begin{enumerate}
	\item the branch point set $\{p_1, \ldots, p_r\}\subset C'$;
	\item generators  $\alpha_1,\ldots,\alpha_{g'},\beta_1,\ldots,\beta_{g'},\gamma_1,\ldots,\gamma_r $ of $ \pi_1(C'\setminus\{p_1, \ldots, p_r\})$, where  each $\gamma_i$ is a simple geometric loop around $p_i$ and 
	 $$\prod_{i=1}^{g'}[\alpha_i,\beta_i]\cdot\gamma_1\cdot\ldots\cdot\gamma_r=1\,;$$
	\item a generating vector $V$ for $H$ of signature $(g';m_1,\ldots ,m_r)$
	with the property that \textit{Hurwitz's formula} holds:
	\begin{equation}\label{Hurwitz}
	2g-2=|H|\bigg(2g'-2+\sum_{i=1}^r\frac{m_i-1}{m_i}\bigg)\,.
	\end{equation}
	\end{enumerate}
	\end{remark}

\begin{remark} Let $C$ be a curve with an action of
a finite group $H$ on it and consider the Galois covering $c\colon C\rightarrow C/H$.
The appropriate orbifold homomorphism $\psi\colon \mathbb T (g';m_1,\ldots, m_r)\rightarrow H$ associated to $V$
is induced by the monodromy of the Galois \'etale $H$-covering
$c^0\colon C^0\rightarrow C'^0$ induced by $c$, where $C'^0$ is the curve obtained from $C'$ by 
removing the branch points $p_1, \ldots, p_r$ of $c$, and $C^0:=c^{-1}(C'^0)$.
 
Let $h_i:= \psi(c_i)\in H$: $h_i$  is called the {\em  local monodromy element} and  generates the stabilizer of a point in $c^{-1}(p_i)$. We also define 
$$\Sigma(V):= \bigcup_{g\in G}\bigcup_{j\in \mathbb Z}\bigcup_{i=1}^r \{ g\cdot h_i^j\cdot g^{-1}\}\,,$$
as the set of stabilizers for the action of $H$ on $C$.
\end{remark}

\begin{remark}\label{disjoint}
 Let $S\rightarrow X:= (C_1 \times C_2)/G$ be a product-quotient surface.
The  two Galois covering $C_i \rightarrow C_i/G=: C'_i$  determine a pair $(V_1,\, V_2)$ of
generating vectors of $G$. 
 The action of $G$  on $C_1 \times C_2$ is free
if and only if $(V_1,V_2)$ is \textit{disjoint}, that is $$\Sigma(V_1)\cap \Sigma(V_2)=\{ 1\}\,.$$

Conversely, a pair $(V_1,V_2)$  of generating vectors of $G$ of  respective signatures 
$T_1=(g'_1; m_1,\ldots, m_r)$ and 
$T_2=(g'_2; n_1,\ldots, n_s)$, determines a family of product-quotient surfaces of dimension $M_1  + r +M_2+s $,
where  $M_i=3g'_i-3$.
\end{remark} 
In particular, the product-quotient surface  can only be  rigid (i.e., has no non-trivial deformations), 
if   $r=s=3$ and $g'_1=g'_2=0$. If the action is non free, then it can happen that the above family has dimension 0, still the surface is not rigid, i.e., it has non trivial deformations which are no longer product-quotient surfaces. 

Instead, if $S$ is isogenous to a product, it follows by the result of the second author below that $S$ is rigid if and only if $r=s=3$ and $g'_1 = g'_2 =0$.

\begin{theorem}[\cite{Cat00}]\label{MT}
a) A projective smooth surface $S$ is isogenous   to a product of two curves of respective genera $g_1, g_2 \geq 2$ ,  if and only if
the following two conditions are satisfied:

\noindent  1) there is an exact sequence
$$
1 \rightarrow \Pi_{g_1} \times \Pi_{g_2} \rightarrow \pi = \pi_1(S)
\rightarrow G \rightarrow 1,
$$
where $G$ is a finite group and where $\Pi_{g_i}$ denotes the fundamental
group of a projective curve of genus $g_i \geq 2$;

\noindent 2) $e(S) (= c_2(S)) = \dfrac{4}{|G|} (g_1-1)(g_2-1)$.

b) Write $ S = (C_1 \times C_2) / G$. Any surface $X$ with the
same topological Euler number and the same fundamental group as $S$
is diffeomorphic to $S$ and is
also isogenous to a product. There is a smooth proper family with connected smooth base manifold $T$, $ p \colon \mathcal{X} \rightarrow T$
having two fibres respectively isomorphic to $X$, and $Y$, where $Y$ is one of the 4 surfaces $ S = (C_1 \times C_2) / G$,
$ S_{+-} :  = (\overline{C_1} \times C_2) / G$, $ \overline{S }= (\overline{C_1}  \times \overline{C_2} ) / G$,
$ S_{-+} :  = (C_1 \times \overline{C_2} ) / G = \overline{ S_{+-}   }$.

c) The corresponding subset of the moduli space of surfaces of general type
$\mathfrak{M}^{top}_S = \mathfrak{M}^{diff}_S$, corresponding to
surfaces orientedly homeomorphic, resp. orientedly diffeomorphic to $S$, is either
irreducible and connected or it contains
two connected components which are exchanged by complex conjugation.

In particular, if $S'$ is orientedly diffeomorphic to $S$, then $S'$ is
deformation equivalent to $S$ or to $\overline{S}$.
\end{theorem}

\begin{definition} Let $S= (C_1 \times C_2)/G$ be a surface isogenous to a product.
$S$ is called a \textit{Beauville surface}   if $r=s=3$ and  $g(C_1/G)=g(C_2/G)=0$.
\end{definition}
A corollary of the above Theorem \ref{MT} (see  \cite{Cat00})  is that a Beauville surface  is rigid,  
it has no nontrivial deformations.

\begin{remark}\label{q=g1+g2}
  Observe  that, by \cite{enr} (cf. also \cite{Ser96}), if $S\rightarrow X:= (C_1 \times C_2)/G$ is a product-quotient
surface, then $q(S) = g(C_1/G) + g(C_2/G)$.		
 Therefore, any Beauville surface is regular.
 \end{remark} 

The concept of  Beauville surfaces was introduced by the second author in 1997, 
and their global rigidity was shown in \cite{Cat00}.
 They  were first systematically studied in \cite{BCG05}, where the connection between  the algebro-geometric background and the group theoretic description was explained in more detail, many new examples were constructed, and a lot of conjectures were presented. 
 Quite a number of these conjectures have been  solved in the meantime and there exists a substantial literature on Beauville surfaces nowadays,
 which is reflected by the several contributions to the present volume.
\section{The fundamental group}

Let $S:=(C_1 \times C_2)/G$  be the minimal realization of  a surface isogenous to a product: then  the surface $C_1\times C_2$ is a Galois \'etale covering space 
of $S$ with group $G$ and we have already observed (in Theorem \ref{MT}) that
 the fundamental group of $S$ sits into an exact sequence
$$1 \longrightarrow \pi_1(C_1)\times \pi_1(C_2) \longrightarrow \pi_1(S)  \longrightarrow G  \longrightarrow 1\,.$$

If one drops the assumption about the freeness of the action of $G$ on the product 
$C_1 \times C_2$, there is no reason that the behaviour of the fundamental group of the quotient should be similar to the 
above situation. Nevertheless, in the unmixed case, it turns out that the fundamental group $\pi_1(X) $
admits a very similar description.

\begin{definition}
We  call the the fundamental group $\Pi_g:=\pi_1(C)$ of a projective curve  of genus $g$ a \textit{(genus $g$) surface group}.
\end{definition}

\begin{theorem}[\cite{BCGP08}]
Let $C_1$ and $C_2$ be projective curves  of genus at least 2 and let $G$ be a finite group acting  on each $C_i$ as a
group of automorphisms and diagonally on the product $C_1 \times C_2$.

Let $S$ be the minimal desingularization of $X:=(C_1 \times C_2)/G$. Then
the fundamental group $\pi_1(X)\cong \pi_1(S)$ has a normal subgroup $\mathcal N$ of finite index 
which is isomorphic to the product of surface groups, i.e., there are integers $h_1, h_2 \geq 0$ such that
$\mathcal N \cong \Pi_{h_1}\times\Pi_{h_2}$.
\end{theorem}

\begin{remark}
1) The previous theorem holds also in dimension $n\geq 3$, see \cite{BCGP08} and \cite{DP10}.\\
2) In the case of surfaces isogenous to a product we have that $h_i$ equals the genus of $C_i$, for $i=1,2$.
\end{remark}

Especially  the  case of infinite
fundamental groups is the one where the above structure theorem  turns out to be extremely helpful
in order  to give an explicit description of these
groups. As we will explain in the sequel, it is easy to get a presentation of the fundamental group of a product-quotient surface. But,
in general,  a presentation of a group does not say much about the group. Even the problem whether it is trivial or not is in general 
an undecidable problem.
Here, since 
we know that $\pi_1$ is a ``surface times surface by finite'' group, we go through  its normal subgroups of finite index 
until such a subgroup $\mathcal N$ appears.

We recall now how to compute the fundamental group of a product-quotient surface $S\rightarrow X= (C_1\times C_2)/G$
starting from the associated algebraic data:
$$\psi_i\colon\mathbb T_i\longrightarrow G\,\qquad i=1,2\,.$$ 

\begin{remark} 
The kernel of $\psi_i$ is isomorphic to the fundamental group $\pi_1(C_i)$, and the action of $\pi_1(C_i)$ on the universal 
cover $u\colon \hat C_i\rightarrow C_i$ extends to a properly discontinuous action of $ \mathbb T_i$. Moreover, $u$ is 
$\psi_i$-equivariant and $\hat C_i/\mathbb T_i\cong C_i/G_i$.
\end{remark}

There are two short exact sequences:
\begin{equation}\label{exseq1}
1\longrightarrow \pi_1(C_i) \longrightarrow \mathbb T_i \stackrel{\psi}{\longrightarrow} G \longrightarrow 1 \,, \quad i=1,2\,.
\end{equation}

\noindent We define the fibre product 
$$\mathbb H :=\mathbb H(G;\psi_1,\psi_2):=\{(x,y)\in \mathbb T_1\times \mathbb T_2 \mid \psi_1(x)=\psi_2(y)\}\,.$$
The exact sequences in  (\ref{exseq1}) induce  then an exact sequence:
$$1\rightarrow \pi_1(C_1)\times \pi_1(C_2) \rightarrow \mathbb H(G;\psi_1,\psi_2) \rightarrow \Delta_G \cong G \rightarrow 1 \,,$$
where $\Delta_G\subset G\times G $ denotes the diagonal subgroup.
\begin{definition}
Let $H$ be a group. Then $\mathrm{Tors}(H)$ is the normal subgroup 
generated by the torsion elements of $H$ (i.e., the elements of finite order in $H$).
\end{definition}

\begin{proposition}[{\cite[Proposition 3.4]{BCGP08}}]
Let $S\rightarrow X= (C_1\times C_2)/G$ be a product-quotient surface and let 
$$\psi_i\colon\mathbb T_i\longrightarrow G\,\quad i=1,2$$ 
be the associated appropriate orbifold homomorphisms.
Then $$\pi_1(S)=\pi_1(X)= \mathbb H/ \mathrm{Tors}(\mathbb H)\,.$$
\end{proposition}

A very special and theoretically easy case is:

\begin{corollary}
Let $S= (C_1\times C_2)/G$ be a surface isogenous to a product.  Then
$\pi_1(S)= \mathbb H$ and $H_1(S, \mathbb Z)= \mathbb H^{ab}$.
\end{corollary}

In the  literature there are several articles on the construction and classification of 
surfaces birational to the quotient of a product of curves; we have already cited most of them.
In some cases the authors provide a computer-script to compute the 
fundamental group (or the first homology group) of the surfaces they study.

Nevertheless, in this note we  take the opportunity to give a simple algorithm  in the special case, 
which includes the case Beauville surfaces, where 
 $S= (C_1 \times C_2)/G$ is a \textit{regular} surface isogenous to a product (equivalently , $g'_1 = g'_2 = 0$).

\

We now recall how to determine whether  two regular surfaces isogenous to a product of unmixed type are deformation equivalent (cf. \cite{BC04, BCG08}).

By Remark \ref{q=g1+g2}, and as observed above, $S$ regular  means that $ C_i/G \cong \mathbb P^1$, $i=1,2$.
We denote by $\mathcal B (G; T_1,T_2)$ the set of disjoint pairs $(V_1,V_2)$ of spherical systems of generators of $G$ of respective signatures
$T_1$ and $T_2$.

Let $V=[h_1,\ldots, h_r]$ be a $r$-tuple of elements of $G$ and $1\leq i	\leq r$. We consider the usual {\em Hurwitz move}
 $\sigma_i(V)$ defined by
$$\sigma_i(V):=[h_1,\ldots, h_{i-1}, h_i h_{i+1} h_i^{-1}, h_i, h_{i+2} ,\ldots, h_r]\,.$$
It is well known  that  $\sigma_1,\ldots, \sigma_r$ generate the braid group on $r$ letters $\mathbf{B}_r$, 
and that $\mathbf{B}_r$ maps spherical system of generators to spherical system of generators, and  preserving the signature.
Also the automorphism group $\mathrm{Aut}(G)$ of $G$ acts on the set of spherical systems of generators of a fixed 
signature by simultaneous application of an automorphism $\varphi$ to the coordinates of a $r$-tuple: 
$\varphi(V):= [\varphi(h_1),\ldots, \varphi(h_r)]$.
We get the following action of $ \mathbf{B}_r\times \mathbf{B}_s \times \mathrm{Aut}(G)$ on $\mathcal{B}(G; T_1,T_2)$:
let $(\gamma_1,\gamma_2,\varphi)\in \mathbf{B}_r\times \mathbf{B}_s \times \mathrm{Aut}(G)$ and
$(V_1,V_2)\in \mathcal{B}(G; T_1,T_2)$, where $T_1$ has length $r$ and $T_2$ has length $s$, we set
$$(\gamma_1,\gamma_2,\varphi)\cdot (V_1,V_2):= (\varphi( \gamma_1(V_1)) ,\varphi( \gamma_2(V_2)))\,.$$
We denote this action by $\mathcal H$.

\begin{theorem}[{cf. \cite[Theorem 5.2]{BCG08}}]
Let $S$ and $S'$ be two regular surfaces, both  isogenous to a product of unmixed type,
and with associate pairs of spherical system of generators $(V_1,V_2)$, respectively  $(V_1',V_2')$.  

Then
$S$ and $S'$ are  deformation-equivalent if and only if the respective groups are isomorphic, $G(S)\cong G(S')$, and 
either $(V_1,V_2)$ and $(V_1',V_2')$ 
 belong to the same  
$\mathcal H$-orbit
or  $(V_1,V_2)$ and $(V_2',V_1')$ do. 
\end{theorem}

\begin{remark}
Given a finite group $G$ and two signatures $T_1$ and $T_2$ of length $r$ and $s$, it is  then work
for a computer to determine $ \mathcal B (G; T_1,T_2) $, its orbits for the  
$\mathcal{H}$-action 
and the first  homolgy group $H_1(S,\mathbb Z)$ of the associate surfaces $S$.
\end{remark}

We have written a MAGMA script (see Appendix \ref{appendix}) which takes as input $G$, $T_1$ and $T_2$,
and returns as output a representative $(V_1,V_2)\in \mathcal{B}(G; T_1,T_2)$ for each 
$\mathcal H$-orbit, together with the first homology group 
$H_1(S,\mathbb Z)$ of the associate surface $S$.
The outline of the script is the following:

\begin{itemize}
\item Step 1: we compute all the spherical systems of generators of $G$ of respective signatures $T_1$ and $T_2$ 
and we collect them in orbits for the action of the braid group $\mathbf{B}_r$ (resp. $\mathbf{B}_s$).

\item Step 2: we discard the pairs of orbits of non-disjoint spherical systems of generators.

\item Step 3: the remaining pairs yield surfaces isogenous to a product.
We consider the action of $Aut(G)$ on them, hence we get the $\mathcal{H}$-orbits in
$ \mathcal B (G; T_1,T_2)  $.
  Observe  that indeed it suffices   to consider only the action of $Out(G)$, since 
 the $Inn(G)$-action was already taken care of  in the first step.

\item Step 4: we run over the outputs of Step 3 and  we compute their first homology group $H_1(S,\mathbb Z)$.

\end{itemize}

\section{Some applications}

Surfaces isogenous to a product of unmixed type with $p_g=0$ have been classified
in \cite{BCG08}.  We run our script for these surfaces and we get:

\begin{theorem}
Let $S = (C_1\times C_2)/G$ be a surface isogenous to a product of unmixed type,
with $p_g(S) = 0$, then $G$ is one of the groups in the Table \ref{tab1} and the signatures are listed in the table.
The number $N$ of components  in the moduli space, their dimension $D$  and the first homology group of $S$ are given in the last three columns.

\begin{table}[!h]
\begin{tabular}{c|c|c|c|c|c|c}
$G$ & $Id(G)$ & $T_1$ & $T_2$ & $N$&$D$&  $H_1(S,\mathbb Z)$\\
\hline
$\mathcal{A}_5$ & $\langle 60,5\rangle$& $[2,5,5]$& $[3,3,3,3]$&1&1&$(\mathbb Z_3)^2\times (\mathbb Z_{15})$\\
$\mathcal{A}_5$ & $\langle 60,5\rangle$& $[5,5,5]$& $[2,2,2,3]$  &1&1&$(\mathbb Z_{10})^2$\\
$\mathcal{A}_5$ & $\langle 60,5\rangle$& $[3,3,5]$& $[2,2,2,2,2]$  &1&2&$(\mathbb Z_2)^3\times \mathbb Z_6$\\
$\mathcal S_ 4 \times \mathbb Z_2$& $\langle 48,48 \rangle$& $[2,4,6]$& $[2,2,2,2,2,2]$  &1&3&
$(\mathbb Z_2)^4\times \mathbb Z_4$\\
G(32) & $\langle 32,27 \rangle$& $[2,2,4,4]$& $[2,2,2,4]$ &1&2&$(\mathbb Z_2)^2\times \mathbb Z_4\times \mathbb Z_8$\\
$(\mathbb Z_5)^2$ & $\langle 25,2\rangle$& $[5,5,5]$& $[5,5,5]$ &1&0&$(\mathbb Z_5)^3$\\
$\mathcal{S}_4$ & $\langle 24,12\rangle$& $[3,4,4]$& $[2,2,2,2,2,2]$ &1&3&$(\mathbb Z_2)^4\times \mathbb Z_8$\\
G(16) & $\langle 16,3\rangle$& $[2,2,4,4]$& $[2,2,4,4]$&1&2&$(\mathbb Z_2)^2\times \mathbb Z_4\times \mathbb Z_8$\\
$D_4\times \mathbb Z_2$ & $\langle 16,11\rangle$& $[2,2,2,4]$& $[2,2,2,2,2,2]$ &1&4&$(\mathbb Z_2)^3\times (\mathbb Z_4)^2$\\
$(\mathbb Z_2)^4 $ & $\langle 16,14\rangle$& $[2,2,2,2,2]$&$[2,2,2,2,2]$  &1&4&$(\mathbb Z_4)^4$\\
$(\mathbb Z_3)^2$ & $\langle 9,2\rangle$& $[3,3,3,3]$& $[3,3,3,3]$&1&2&$(\mathbb Z_3)^5$\\
$(\mathbb Z_2)^3$ & $\langle 8,5\rangle$& $[2,2,2,2,2]$& $[2,2,2,2,2,2]$  &1&5&$(\mathbb Z_2)^4\times (\mathbb Z_4)^2$\\
\end{tabular}	
\caption{ }\label{tab1}
\end{table}
\end{theorem}

\begin{remark} (1) In a  previous paper \cite{BC04} (containing the classification of 
surfaces isogenous to a product of unmixed type with $p_g=q=0$ and $G$ abelian) there is a mistake, 
an erroneous statement about commutators in a product. The error resulted
into finding only a proper quotient of the actual first homology group.

(2) The correct calculation of the first homology group was done more
than one year ago by the first author, who used a MAGMA script to
find the correct answer.

(3) Another correction to do to \cite{BC04} is that the authors forgot the
possibility of swapping factors, hence there is only one isomorphism
class of surfaces in the case $G = (\mathbb Z_5)^2$.

(4) Our computer calculations confirm the result of the paper
\cite{Shab13}, dedicated to the calculation of the first homology groups of surfaces isogenous to 
a product of unmixed type with $p_g=q=0$ and $G$ abelian.
\end{remark}

In \cite{Glei13}, the author classifies the regular surfaces isogenuos to a product of unmixed type with $\chi (\mathcal O)=2$; 
in particular he classifies the unmixed Beauville surfaces with $p_g=1$.
Applying our program we get the following:
\begin{theorem}[cf. {\cite{Glei13}}]
Let $S = (C_1\times C_2)/G$ be a Beauville surface of unmixed type,
with $p_g(S) = 1$. Then $G$ is one of the groups in  Table \ref{tab2} and the signatures are listed in the table.
The number  $N$ of components in the moduli space  and the first homology group of $S$ are given in the last two columns.
\begin{table}[!h]
\begin{tabular}{c|c|c|c|c|c}
$G$ & $Id(G)$ & $T_1$ & $T_2$ & $N$&  $H_1(S,\mathbb Z)$\\
\hline
$PSL(2,7)\times \mathbb Z_2$ & $\langle 336, 209\rangle$& $[2,3,14]$& $[4,4,4]$&2&$(\mathbb Z_4)^2$\\
$PSL(2,7)$ & $\langle 168,42\rangle$& $[7,7,7]$& $[3,3,4]$  &2&$\mathbb Z_7\times \mathbb Z_{21}$\\
$PSL(2,7)$ & $\langle 168,42\rangle$& $[3,3,7]$& $[4,4,4]$  &2&$\mathbb Z_4\times \mathbb Z_{12}$\\
$G(128,36)$& $\langle 128,36 \rangle$& $[4,4,4]$& $[4,4,4]$  &2&$(\mathbb Z_2)^3\times (\mathbb Z_4)^2$\\
\end{tabular}	
\caption{ }\label{tab2}
\end{table}
\end{theorem}
Let $G=\mathrm{PSL}(2,q)$, where $q$ and is a prime power and $q\leq9$.
For $q\leq 5$, no group has a disjoint pair of spherical generators (see \cite{BCG05}), while in the other cases
the ``admissible'' pairs of signatures are classified (see \cite{gar10}).
The outputs of our script are collected in Table \ref{tab3}. We use the same notation of the previous tables and
we add a column  reporting the Euler characteristic $\chi(\mathcal O_S)$.

\begin{table}[!h]
\begin{tabular}{c|c|c|c|c|c|c}
$G$ & $Id(G)$ & $T_1$ & $T_2$ & $N$&  $H_1(S,\mathbb Z)$& $\chi(\mathcal{O}_S)$\\
\hline
\multirow{ 5 }{*}{$PSL(2,7)$} & \multirow{ 5 }{*}{$\langle 168,42\rangle $} 
    &$[ 3,3,4]$ & $[7,7,7]$ & 2& $\mathbb Z_7 \times \mathbb Z_{21}$&2\\
 &&$[ 3,4,4]$ & $[7,7,7]$ & 1& $\mathbb Z_7 \times \mathbb Z_{28}$&4\\
 &&$[ 4,4,4]$ & $[7,7,7]$ & 2& $(\mathbb Z_{28})^2$&6\\
 &&$[ 4,4,4]$ & $[3,7,7]$ & 4& $\mathbb Z_4 \times \mathbb Z_{28}$&4\\
 &&$[ 4,4,4]$ & $[3,3,7]$ & 2& $\mathbb Z_4 \times \mathbb Z_{12}$&2\\
\hline
\multirow{ 8 }{*}{$PSL(2,8)$} & \multirow{ 8 }{*}{$\langle 504,156 \rangle $} 
    & $[2,7,7]$ & $[3,3,9]$ & 3 &$\mathbb Z_3 \times \mathbb Z_{21}$&6\\  
&  & $[2,7,7]$ & $[3,9,9]$ & 3 &$\mathbb Z_3 \times \mathbb Z_{63}$&12\\
&  & $[2,7,7]$ & $[9,9,9]$ & 7& $\mathbb Z_9 \times \mathbb Z_{63}$&18\\
&  & $[7,7,7]$ & $[9,9,9]$ & 14 &$(\mathbb Z_{63})^2$ &48\\
&  & $[7,7,7]$ & $[2,9,9]$ & 6& $\mathbb Z_7 \times \mathbb Z_{63}$&20\\
&  & $[7,7,7]$ & $[2,3,9]$ & 6&$\mathbb Z_7 \times \mathbb Z_{21}$&4\\
&  & $[7,7,7]$ & $[3,9,9]$ & 6&$\mathbb Z_{21} \times \mathbb Z_{63}$&32\\
&  & $[7,7,7]$ & $[3,3,9]$ & 6&$(\mathbb Z_{21})^2$&16\\
\hline
 \multirow{ 6 }{*}{$PSL(2,9)\cong \mathcal A_6$} & \multirow{ 6 }{*}{$\langle 360,118\rangle $} 
&$[ 3,3,5]$ & $[4,4,4]$ & 1& $(\mathbb Z_{12})^2 $&3\\
 &  &$[ 3,5,5]$ &$ [4,4,4]$ & 1& $\mathbb Z_4 \times \mathbb Z_{60}$&6\\
  &&\multirow{ 2 }{*}{$[ 5,5,5]$} &\multirow{ 2 }{*}{ $[4,4,4]$ }& 2& $\mathbb Z_{20}\times \mathbb Z_{60}$&\multirow{2}{*}{9}\\
 &&&&2& $(\mathbb Z_{20})^2$&\\	
&&\multirow{ 2 }{*}{$[ 5,5,5]$}& \multirow{ 2 }{*}{$[3,3,5]$} & 2& $\mathbb Z_5 \times \mathbb Z_{15}$&\multirow{2}{*}{3}\\
 &&& & 2& $(\mathbb Z_{15})^2$&\\
\end{tabular}	
\caption{ }\label{tab3}
\end{table}

In Table \ref{tab4}, we collect the outputs in some other easy cases: $G\in\{ \mathcal{S}_5, \mathcal{S}_6,
(\mathbb{Z}_7)^2\}$.
One can prove (using a simple MAGMA script), that if $X=(C_1\times C_2)/G$ is a Beauville surface with 
 $G\in\{ \mathcal{S}_5, \mathcal{S}_6,(\mathbb{Z}_7)^2\}$, then its signature $(T_1,T_2)$ is in Table \ref{tab4}.
 
\begin{table}[!h]
\begin{tabular}{c|c|c|c|c|c|c}
$G$ & $Id(G)$ & $T_1$ & $T_2$ & $N$&  $H_1(S,\mathbb Z)$& $\chi(\mathcal{O}_S)$\\
\hline
$\mathcal{S}_5$&$\langle120,34 \rangle$ &$[4,4,5]$&$[3,6,6]$ &1 &$\mathbb Z_{3} \times \mathbb Z_{24} $&3\\
\hline
\multirow{ 5 }{*}{$\mathcal{S}_6$}& \multirow{ 5 }{*}{$\langle 720 ,763\rangle $} 
																	&$[5,6,6 ]$, & $[4,6,6]$ & 8 &$(\mathbb Z_6)^2 $ &35\\
                  &                                  &$[5,6,6 ]$, & $[4,4,6]$ &16 &$ \mathbb Z_{2} \times \mathbb Z_{24}$& 28\\
                  &                                  &$[5,6,6 ]$, & $[4,4,4]$ &  8 &$\mathbb Z_{4} \times \mathbb Z_{12} $ &21\\
                  &                                  &$[3,6,6 ]$, & $[4,4,4]$ & 5 &$(\mathbb Z_{12})^2 $ &15\\
                  &                                 	&$[2,5,6 ]$, & $[4,4,4]$ & 1 &$(\mathbb Z_4)^2 $ &6\\
\hline
$(\mathbb{Z}_7)^2$&$\langle 49,2 \rangle$ &$[7,7,7]$&$[7,7,7]$ &7 &$(\mathbb Z_7)^3$&4\\
\end{tabular}	
\caption{}\label{tab4}
\end{table}

\appendix

\section{ The script}\label{appendix}
{\small \begin{verbatim}

// TuplesOfGivenOrder creates a sequence of the same length as the input
// sequence type, whose entries are subsets of the group in the input,
// and precisely the subsets of elements of order the corresponding
// entry of type

TuplesOfGivenOrders:=function(G,type)
SEQ:=[];
for i in [1..#type-1] do
  EL:={g: g in G| Order(g) eq type[i]};
  if IsEmpty(EL) then return [{}];
  else Append(~SEQ,EL);
end if; end for;
return SEQ;
end function;

// This transforms a tuple into a sequence

TupleToSeq:=func<tuple|[x: x in tuple]>;

// Now we create all sets of spherical generators of a group of the
// prescribed signature.

VectGens:=function(G, type)
Vect:={}; SetCands:=TuplesOfGivenOrders(G,type);
for cand in CartesianProduct(SetCands) do
  if Order(&*cand) eq type[#type] then
    if #sub<G|TupleToSeq(cand)> eq #G then 
      Include(~Vect, Append(TupleToSeq(cand),(&*cand)^-1));
end if; end if; end for; ;
return Vect;
end function;

// HurwitzOrbit, starting from a sequence seq of elements of a group,
// creates all sequences of elements which are equivalent to the given one
// for the equivalence relation generated by the Hurwitz moves
// and returns (to spare memory) only the ones whose entries have 
// orders disposed as the ones in seq

HurwitzMove:= func<seq,j|Insert(Remove(seq,j),j+1,seq[j]^seq[j+1])>;

HurwitzOrbit:=function(seq)
orb:={ }; shortorb:={  }; Trash:={ seq };
repeat ExtractRep(~Trash,~gens); Include(~orb, gens);
  for k in [1..#seq-1] do 
    newgens:=HurwitzMove(gens,k);
    if newgens notin orb then Include(~Trash, newgens);
end if; end for; until IsEmpty(Trash);
for gens in orb do  test:=true;
  for k in [1..#seq] do 
    if not Order(gens[k]) eq Order(seq[k]) then 
      test:=false; break k;
  end if; end for;
  if test then Include(~shortorb, gens); end if;
end for;
return shortorb;
end function;

// OrbitsVectGens creates a sequence 
// whose entries are subsets of the group in the input,
// and precisely the subsets  are the orbits under the
// Hurwitz moves

OrbitsVectGens:=function(G, type)        
Orbits:=[];   Vects:=VectGens(G, type);
while not IsEmpty(Vects) do
  v:=Rep(Vects);  orb:=HurwitzOrbit(v);
  Append(~Orbits, orb); Vects:=Vects diff orb;
end while;
return Orbits;
end function;
 
// Disjoint checks if two spherical systems of generators are disjoint 

Stab:= function(seq, G)
M:= {} ;
for i in [1..#seq] do g:=seq[i]; 
  for n in [1 .. (Order(g)-1) ] do
    M := M join Conjugates(G,g^n) ;
end for; end for;
return M;
end function ;

Disjoint:=func<G,seq1, seq2 | IsEmpty( Stab(seq1,G) meet Stab(seq2,G))>;

// Homology computes the first homology group
// of the surface associated to the disjoint pair
// of spherical systems of generators (seq1, seq2) of G

Poly:=function(seq, gr) 
F:=FreeGroup(#seq);  Rel:={}; Q:=Id(F);
for i in {1..#seq} do 
  Q:=Q*F.(i); Include(~Rel,F.(i)^(Order(seq[i])));   
end for;
Include(~Rel,Q); P:=quo<F|Rel>; 	
return P, hom<P->gr|seq>;
end function;

Homology:=function(G, seq1,  seq2)
T1,f1:=Poly(seq1,G);   T2,f2:=Poly(seq2,G);
T1xT2:=DirectProduct(T1,T2);
 GxG,inG:=DirectProduct(G,G);
if Category(G) eq GrpPC then 
  n:=NPCgens(G); 
  else n:=NumberOfGenerators(G); 
end if;
Diag:= hom<G->GxG| [inG[1](G.j)*inG[2](G.j): j in [1..n]]>(G);
f:=hom<T1xT2->GxG| inG[1](seq1) cat inG[2](seq2)>;
Pi1:=Rewrite(T1xT2,Diag@@f); 
return AbelianQuotient(Pi1);
end function;

// Surfaces is the main function.
// It takes as input a group G and two signatures.
// It calls the previous function and moreover identifies
// distinct pairs of Hurwitz's orbits of spherical systems of
// generators under the action of Aut(G). 
// Note that we consider only the action of Out(G), since 
// the Inn(G)-action was already taken care of in HurwitzOrbit.

Surfaces:=function(G, type1,type2) 
R:={}; Aut:=AutomorphismGroup(G);
F,q:=FPGroup(Aut);     O1,p:=OuterFPGroup(Aut);
Out,k:=PermutationGroup(O1);
V1:=OrbitsVectGens(G, type1); 
if type1 eq type2 then V2:=V1;
  else V2:=OrbitsVectGens(G,type2); 
end if;
W:=Set(CartesianProduct({1..#V1},{1..#V2}));
for pair  in W do 
  v1:= Rep(V1[pair[1]]); v2:= Rep(V2[pair[2]]);
  if not Disjoint(G, v1,v2) then Exclude(~W, pair);
 end if;
end for;
while not IsEmpty(W) do 
  pair:=Rep(W);    
  v1:= Rep(V1[pair[1]]); v2:= Rep(V2[pair[2]]);Include(~R,[v1,v2]);
  for phi in Out do 
    if exists(x){y: y in W | (q(phi@@(p*k)))(v1) in V1[y[1]] and 
   	 (q(phi@@(p*k)))(v2) in V2[y[2]]} then Exclude(~W, x); 
      if type1 eq type2 then Exclude(~W, <x[2],x[1]>); 
    end if; end if;
    if IsEmpty(W) then break phi; end if;
end for; end while;
printf "Number of components: %o\n", #R;
for r in R do 
  printf "Spherical generators:\n%o\n%o\n", r[1],r[2];
  printf "Homology:\n%o\n", Homology(G,r[1],r[2]);
end for;
return {};
end function;
\end{verbatim}
}

\

\noindent{\bf Author's Adress:}\\
Ingrid Bauer, Fabrizio Catanese, Davide Frapporti\\
Lehrstuhl Mathematik VIII\\
Mathematisches Institut der Universit\"at Bayreuth, NW II\\
Universit\"atsstr. 30; D-95447 Bayreuth, Germany

\end{document}